\documentclass[a4paper,11pt]{article}

\usepackage[centertags]{amsmath}
\usepackage[nosort]{cite}
\usepackage[psamsfonts]{amsfonts} 
\usepackage{amsthm,enumerate,amssymb}

\newcommand{\dd}{\mathrm{d}}

\newcommand{\bigtimes}{\mathop{\mbox{\Large $\times$}}} 
\newcommand{\QED}{\hfill{$\square$} \\ \medskip}
\newcommand{\tr}{\operatorname{tr}}
\newcommand{\bbbone}{{\mathchoice {\rm 1\mskip-4mu l} {\rm 1\mskip-4mu l}
{\rm 1\mskip-4.5mu l} {\rm 1\mskip-5mu l}}}
\newcommand{\is}{\mathcal{H}}
\newcommand{\LZ}{\mathcal{L}^\infty_{\mathbb{Z}}}
\newcommand{\Lm}{\mathcal{L}^\infty_{[-m,m]}}
\newcommand{\DZ}{\mathcal{D}_{\mathbb{Z}}}

\theoremstyle{plain}
\newtheorem{theorem}{Theorem}

\begin{document}

\title{A Universal Dilation of Discrete Markov Evolutions}

\author{M.~Gregoratti\thanks{E-mail:
matteo.gregoratti@polimi.it}\\ \\
{\footnotesize \textsl{Dipartimento di Matematica ``F.Brioschi'', Politecnico di
Milano}} \\
{\footnotesize \textsl{Piazza Leonardo da Vinci 32, I-20133 Milano, Italy}}\\ \\
{\small Quaderno di Dipartimento n.\textbf{649/P}}}

\date{December 15, 2005}

\maketitle

\begin{abstract}
Given a finite state space $E$, we build a universal dilation for all possible discrete time Markov 
chains on $E$, homogeneous or not: we introduce a second system (an ``environment'') and a 
deterministic invertible time-homogeneous global evolution of the system $E$ with this environment
such that
any Markov evolution of $E$ can be realized by a proper choice of the initial (random) state of the
environment, which therefore determines the transition probabilities of the system. 
We also compare this dilation with the quantum dilations of a Quantum Dynamical Semigroup:
given a Classical Markov Semigroup, we show that it can be extended to a Quantum Dynamical 
Semigroup for which we can find a quantum dilation to a group of $*$-automorphisms
admitting an invariant abelian subalgebra
where this quantum dilation gives just our classical dilation.

\bigskip


\bigskip

\end{abstract}

\section{INTRODUCTION}

A well known theorem due to Birkhoff shows that any doubly stochastic matrix is a convex combination
of permutation matrices. These describe with matrix terminology one-step evolutions which are both
deterministic and invertible. Analogously, a theorem due to A.~S.~Davis \cite{D61}
shows that any stochastic 
matrix, doubly stochastic or not, is a convex combination of deterministic matrices. These describe 
with matrix terminology one-step evolutions which are deterministic but not necessarily invertible. As
a consequence, A.~S.~Davis shows that each (finite state) Markov chain can be realized as an automaton
with random inputs, thus establishing a connection between Markov chains and automata theory.

The relationship between Markov chains and deterministic dynamics can be further analyzed looking for
realizations of Markov evolutions as deterministic invertible time-homogeneous evolutions of the
system coupled with a second system.  
Such representations are usually called dilations and they are
important to show that a Markov evolution of a given system is compatible with a deterministic 
invertible homogeneous model for the evolution of a bigger system. Thus the stochastic features of the
system evolution are given a dynamical explanation. This is theoretically relevant if Markov chains
are applied to phenomena, like physical phenomena for example, for which an underlaying theory
postulates deterministic invertible homogeneous evolutions in absence of noise and external
disturbances. For these phenomena the second system introduced by the dilation models the surrounding
world, the environment, which is the source of the noise.

Dilations are indeed typical in Quantum Probability, where the cultural influence of Quantum Physics
always emphasized the importance that a stochastic evolution of a (non isolated) quantum
system does not contradict the axioms of Quantum Mechanics, i.e.\ that it
can arise form a unitary homogeneous evolution of a (isolated) bigger quantum system,
consisting of the given system and its environment.

Here we are interested in the classical analogue of quantum dilations. 
More precisely, in a completely classical
framework, we consider a system with finite state space $E$, undergoing a discrete 
time evolution given by a Markov chain, not necessarily homogeneous. 
Then, in order to get a dilation of this chain, a second system is introduced, called environment, 
with its state space $(\Gamma, \mathcal{G})$, a measurable space, together with a global
invertible one-step evolution $\alpha:E\times\Gamma\to E\times\Gamma$. Thus, if $(i,\gamma)$ is the 
state of the
compound system at time 0, then $\alpha^t(i,\gamma)$ is its state at time $t\in\mathbb{N}$,
where hence $\alpha^t$ gives a deterministic invertible homogeneous global evolution. If,
nevertheless, the state of the environment is never observed and if initially it is randomly
distributed with some law $Q$ on $(\Gamma,\mathcal{G})$, the evolution of the observed system turns 
out to be stochastic and, if $\Gamma$, $\mathcal{G}$,
$\alpha$ and $Q$ are properly built, it is given by the original Markov chain.
In this case, we say that 
$(\Gamma,\mathcal{G}, \alpha, Q)$ is a dilation of the Markov evolution on $E$.

In this paper, given only the state space $E$ (arbitrary but finite), we build a \emph{universal} 
dilation
$(\Gamma,\mathcal{G},\alpha,\{Q\})$, where $\{Q\}$ is an entire family of distributions 
which can produce any Markov chain on $E$: 
every Markov chain, homogeneous or not, can be dilated by taking always the same, universal, 
$(\Gamma,\mathcal{G},\alpha)$ and by choosing every time the proper distribution $Q$
for the initial state of the environment.
Moreover, our construction, which relays on A.~S.~Davis decomposition of stochastic matrices, allows
to interpret each Markov chain, not only as the stochastic
dynamics resulting from the coupling with an
environment, but at the same time also as an automaton with random inputs, which are now dynamically
provided by the environment itself.

An unpublished result by K\"ummerer \cite{K82} provides a dilation of a Markov chain
with a construction similar to ours, but it does not exhibit the same universality because
the interaction $\alpha$ depends on the chain under consideration.
Our aim is similar also to the aim of Lewis and Maassen \cite{LM84} when they consider classical
mechanics in continuous time and, taken a linear Hamiltonian system modelling a particle and its
environment, they describe how Gibbs states of the whole system lead to stationary Gaussian stochastic
processes for the observables pertaining to the particle under consideration. However, we do not look
for good global states, but for good states $Q$ of the environment alone which lead to Markov
evolutions of the system $E$, our particle, for every independent choice of its initial state.

The paper also shows the relationship between a universal dilation 
$(\Gamma,\mathcal{G},\alpha,\{Q\})$ and the quantum constructions which inspired it. In Quantum 
Probability \cite{G,K82,K88,M,P,P89} the starting object is a Quantum Dynamical Semigroup (QDS) which
describes the evolution of a non isolated system. It is the quantum analogue of a Classical Markov 
Semigroup (CMS). Then a dilation realizes the QDS as a group of $*$-automorphisms describing an 
invertible homogeneous evolution of the system coupled with an environment or, equivalently,
as the expectation of a Quantum Stochastic Flow, which is a quantum generalization of a Markov
stochastic process. We shall prove that our dilations $(\Gamma,\mathcal{G},\alpha,Q)$ are 
classical restrictions of quantum dilations: every CMS on $E$ admits an extension to a QDS for 
which we can find a quantum dilation which is itself an extension of the dilation of the CMS.
However, we shall not embed a universal dilation
$(\Gamma,\mathcal{G},\alpha,\{Q\})$ in the quantum world, as quantum dilations do not exhibit the same
universality and they strictly depend on the QDS under consideration, so that it is not enough
to change the environment state to get another QDS.

The paper is divided into two parts. The first part is completely classical and it consists of
Section~2, where we set up the framework, and of Section~3, where we give the definition of 
universal dilation of Markov evolutions on
$E$ and we show that there exists always one by a constructive proof. Then the second part, Section~4,
is devoted to the set up of the quantum framework and the analysis of the relation
between our (classical) universal dilation of a CMS
and the (quantum) dilation
of a QDS.

\section{PRELIMINARIES}
We consider a system with finite state space $E=\{1,\ldots,N\}$, fixed for the whole paper.
We denote by $P=(p_{ij})_{i,j\in E}$ a stochastic matrix on $E$, so that $p_{ij}\geq0$ and $\sum_j
p_{ij}=1$ for every $i$, and we identify the elements of the (complex abelian $*$-) 
algebra 
$\mathcal{L}^\infty(E)$, the system random variables
$f:E\to\mathbb{C}$, with the column vectors in $\mathbb{C}^N$, so that every stochastic matrix $P$ on
$E$ defines an operator on $\mathcal{L}^\infty(E)$,
\begin{equation*}
    \big(Pf\big) (i)=\sum_{j\in E}p_{ij}\,f(j),
\end{equation*}
describing a one-step probabilistic evolution. Taken a sequence of stochastic matrices
$(P(t))_{t\in\mathbb{N}}$, the evolution of a system random variable $f$
from time 0 to time $t\geq1$ is therefore given by
\begin{equation}\label{cobsev0}
    f \mapsto f_t := P(1) \cdots P(t) \,f, \qquad \forall f\in \mathcal{L}^\infty(E).
\end{equation}
If the sequence is constant, $P(t)=P$ for every $t$, then
the evolution is homogeneous and it is described by the Classical Markov
Semigroup $(P^t)_{t\geq0}$.

We denote by $D$ a deterministic matrix on $E$, a stochastic matrix with
a 1 in each row. Every $D$ describes with matrix terminology a deterministic (not necessarily
invertible) evolution $\beta$, where
\begin{equation}\label{detmatrix}
    D=(d_{ij})_{i,j}, \qquad \qquad \beta:E\to E, \qquad \qquad d_{ij}=\delta_{\beta(i),j},
\end{equation}
so that $Df=f\circ\beta$.
The invertible (bijective) maps on $E$ correspond to the special cases of permutation matrices.
The deterministic matrices are just the extreme points of the convex set of stochastic matrices
and every $P$ is a convex combination of deterministic matrices,
\begin{equation}\label{dcc}
    P=\sum_{\ell\in L}q_\ell \, D_\ell, \qquad \qquad 
    q_\ell\geq0, \quad \sum_{\ell\in L} q_\ell=1.
\end{equation}
One can find such a decomposition with $N^N$ terms, with the set $L$ labelling all possible
deterministic matrices and weighing each $D_\ell$ with $q_\ell=p_{1\beta_\ell(1)}\cdots 
p_{N\beta_\ell(N)}$. Let us remark that the decomposition \eqref{dcc} is not unique and that, for 
any given
$P$, no more than $N^2-N+1$ terms are needed \cite{D61}. Anyway, since we are not going to fix $P$, we
shall employ \eqref{dcc} in the described form: a sum of $N^N$ terms which can produce any $P$ simply
by changing the weights $(q_\ell)_{\ell\in L}$.

For every sequence of stochastic matrices $(P(t))_{t\in\mathbb{N}}$, there exists a Markov chain\\
$\big(\Omega,\mathcal{F},(\mathcal{F}_t)_{t\geq0},(X_t)_{t\geq0},(\mathbb{P}_k)_{k\in E}\big)$ 
with transition probabilities given by $P(t)$, i.e. a discrete time
stochastic process of random variables
$X_t:\Omega\to E$, adapted to a filtration $(\mathcal{F}_t)_{t\geq0}$, and a family of probability
measures $\mathbb{P}_k$, $k\in E$, such that the starting distribution of the process depends on $k$,
$X_0$ has Dirac distribution $\delta_k$ under $\mathbb{P}_k$, but the process always satisfies the 
Markov property with transition 
matrices $P(t)$:
\begin{equation*}
    \mathbb{P}_k(X_{t+1}=j|\mathcal{F}_t) = \mathbb{P}_k(X_{t+1}=j|X_t)
    = p_{X_tj}(t+1), \qquad \forall k,j,t.
\end{equation*}
Chosen a starting distribution $\delta_k$, the Markov process law is then uniquely determined by 
the transition matrices $P(t)$.
Moreover, a system random variable $f\in \mathcal{L}^\infty(E)$ has now a stochastic evolution
described by the $*$-unital homomorphism
\begin{equation}\label{stsystev}
    j_t:\mathcal{L}^\infty(E)\to\mathcal{L}^\infty(\mathcal{F}_t), \qquad f\mapsto
    j_t[f]:=f(X_t), \qquad t\geq0,
\end{equation}
and the evolution \eqref{cobsev0} admits the representation
\begin{equation}\label{cobsev1}
    f_t(k) = \Big(P(1) \cdots P(t) \,f\Big) (k) = \mathbb{E}_k\big[f(X_t)\big], 
    \qquad \forall f\in \mathcal{L}^\infty(E),
\end{equation}
and we say that the system has a Markov evolution.

Let us briefly show a realization of
$\big(\Omega,\mathcal{F},(\mathcal{F}_t)_{t\geq0},(X_t)_{t\geq0},(\mathbb{P}_k)_{k\in E}\big)$ 
by means of the decomposition \eqref{dcc}. If $P(t)=\sum_{\ell\in L}q_\ell(t) \, D_\ell$, then one can
take
\begin{gather*}
    \Omega=E\times L^{\mathbb{N}},\;\;\;
    \omega=(i,\ell_1,\ell_2,\ldots),\;\;\;
    \mathcal{F}=\text{cylindric $\sigma$-field},\;\;\;
    \mathcal{F}_t=\sigma\big(\{i,\ell_1,\ldots,\ell_t\};i\in E,\ell_n\in L \big),\\
    X_0(\omega)=i,\qquad
    X_t(\omega)=\beta_{\ell_t}\circ\cdots\circ\beta_{\ell_1}(i), \qquad
    \mathbb{P}_k\big(\{i,\ell_1,\ldots,\ell_t\}\big)=\delta_{i,k}\prod_{s=1}^t q_{\ell_s}(s),
\end{gather*}
where $\{i,\ell_1,\ldots,\ell_t\}$ is the event 
$\{\omega=(i,\ell_1,\ell_2,\ldots,\ell_t,\ell_{t+1},\ldots),\; \ell_s\in L, s>t\}$.
In this way the Markov chain associated to $(P(t))_{t\in\mathbb{N}}$ is represented as an 
automaton with independent
random inputs: at every step the set of all possible mappings $\beta_\ell:E\to E$ is available and the
system evolution is determined by the value of the input parameter $\ell$ which is selected randomly
according to $q(t)$ and independently of the previous steps. One could also explicitly introduce the
random input parameters
\begin{equation*}
    Y_t:\Omega\to L, \qquad Y_t(\omega)=\ell_t, \qquad t\in\mathbb{N},
\end{equation*}
which are clearly independent under every $\mathbb{P}_k$. Then 
\begin{equation*}
    \mathcal{F}_t=\sigma(X_0, Y_s ,1\leq s\leq t), \qquad X_t=\beta_{Y_t}(X_{t-1}).
\end{equation*}
If the sequence $P(t)$ is
constant, then the Markov chain is homogeneous and the input parameters can be chosen identically
distributed.
To get a dilation, we shall introduce a bigger $\Omega$, namely $E\times G^\mathbb{Z}$, with a $G$
bigger than $L$ in order to define an invertible one-shot coupling $\varphi$, and with $\mathbb{Z}$ 
instead of $\mathbb{N}$ in order to get an invertible dynamic with group properties $\alpha^t$.

In the sequel, given a function $f$ on a domain $E$, we shall always denote with the same symbol $f$
also its trivial extension  to a domain $E\times\Gamma$, $f(i,\gamma)=f(i)$.

\section{UNIVERSAL DILATION OF DISCRETE MARKOV EVOLUTIONS ON A FINITE STATE SPACE}

Let us denote simply by $\{P\}$ a sequence $(P(t))_{t\in\mathbb{N}}$ of stochastic matrices on the
state space $E$ and let us denote by $\mathcal{P}$ the set of such sequences.
We call \textbf{universal dilation} of the Markov evolutions on $E$ a term
\begin{equation}\label{ud}
    \Big(\Omega,\mathcal{F},(\mathcal{F}_t)_{t\geq0},\big(Z_t=(X_t,\Upsilon_t)\big)_{t\geq0},
    (\mathbb{P}_{k,\{P\}})_{k\in E,\{P\}\in\mathcal{P}}\Big),
\end{equation}
where
\begin{itemize}
\item $\big(\Omega,\mathcal{F},(\mathcal{F}_t)_{t\geq0},\mathbb{P}_{k,\{P\}}\big)$
is a filtered probability space for every $k$ and $\{P\}$,
\item every $Z_t=(X_t,\Upsilon_t)$ is a random variable with values in $E\times\Gamma$, being
$(\Gamma,\mathcal{G})$ a fixed measurable space,
\item for every choice of $k$ and $\{P\}$,
$X_0$ has distribution $\delta_k$, independently of $\{P\}$, while $\Upsilon_0$ has distribution
independent of $k$, say $Q_{\{P\}}$,
\item for every choice of $\{P\}$, the term
$\big(\Omega, \mathcal{F}, (\mathcal{F}_t)_{t\geq0}, (X_t)_{t\geq0}, (\mathbb{P}_{k,\{P\}})_{k\in
E}\big)$ is a Markov chain with transition matrices $\{P\}$,
\item there exists an invertible measurable map $\alpha:E\times\Gamma\to E\times\Gamma$ such that
$Z_t=\alpha^t(Z_0)$ for every $t\geq0$.
\end{itemize}

Thus, besides the system $E$, a second system is introduced, an environment with state space
$(\Gamma,\mathcal{G})$. Their states $X_t$ and $\Upsilon_t$ are asked to be random variables on a 
same measurable
space $(\Omega,\mathcal{F})$ such that the global state $Z_t=(X_t,\Upsilon_t)$ undergoes a 
deterministic invertible homogeneous evolution $\alpha^t$. Therefore all the $X_t$ and $\Upsilon_t$ 
are determined by $Z_0$, so that $X_t$ and $\Upsilon_t$ are measurable with respect to
$\sigma(Z_0)=\sigma(X_0,\Upsilon_0)\subseteq\mathcal{F}$ and, depending on the probability 
chosen on $\mathcal{F}$, they are deterministic if and only if $Z_0$ is. The probabilities 
$\mathbb{P}_{k,\{P\}}$ actually fix the value of $X_0$ but usually not that of
$\Upsilon_0$. They depend on two indexes, $k$ which fixes the value of $X_0$, and $\{P\}$
which fixes the distribution of $\Upsilon_0$ and the resulting properties of 
the $X_t$.
The space $(\Omega, \mathcal{F})$ is also endowed with a filtration $\mathcal{F}_t$. 
Note that only the $X_t$ are assumed to be adapted to $\mathcal{F}_t$ so that, in particular,
$\Upsilon_0$ does not have to be $\mathcal{F}_0$-measurable. Therefore
the $X_t$ are not trivially $\mathcal{F}_0$-measurable, even if their values 
are completely determined by the values of $X_0$ and 
$\Upsilon_0$, and,
neglecting the environment, each $\big(\Omega, \mathcal{F}, (\mathcal{F}_t)_{t\geq0}, (X_t)_{t\geq0}, 
\mathbb{P}_{k,\{P\}}\big)$ can be a process with a proper stochastic evolution. What we ask
is that $\big(\Omega, \mathcal{F}, (\mathcal{F}_t)_{t\geq0}, (X_t)_{t\geq0}, 
\mathbb{P}_{k,\{P\}}\big)$ actually is a Markov chain starting from $k$ with transition matrices 
$\{P\}$. At the same time however, such a dilation shows that this Markov evolution of the system 
$E$ is compatible 
with a deterministic, invertible and homogeneous model for the evolution of $E$ coupled with an
environment $\Gamma$. In particular, as $X_0=k$, the whole stochasticity of the process is due only 
to the randomness of the unobserved environment initial state $\Upsilon_0$. 

We call universal such a 
dilation because the same
$\Big(\Omega,\mathcal{F},(\mathcal{F}_t)_{t\geq0},\big(Z_t=(X_t,\Upsilon_t)\big)_{t\geq0}\Big)$ allows
to represent all the Markov evolutions on $E$ with the change of the probability
$\mathbb{P}_{k,\{P\}}$ alone. Therefore, both the environment state space $(\Gamma,\mathcal{G})$ and 
the global evolution
$\alpha$ do not depend on the Markov chain to be dilated, but only on the state space $E$.

A universal dilation gives another interpretation of every evolution \eqref{cobsev0}, 
compatible with the representation \eqref{cobsev1}:
\begin{equation}\label{cobsev2}
    f_t(k) = \Big(P(1) \cdots P(t) \,f\Big)(k) = \mathbb{E}_{k,\{P\}}\big[f(Z_t)\big] = 
    \mathbb{E}_{k,\{P\}}\Big[f\big(\alpha^t(k,\Upsilon_0)\big)\Big], \qquad \forall f\in
    \mathcal{L}^\infty(E).
\end{equation}
The stochastic evolution of a system variable $f\in \mathcal{L}^\infty(E)$ is now
described by $f\mapsto j_t[f]:=f(X_t)=f(Z_t)=f\circ\alpha^t(Z_0)$. 
And we can also consider global
random variables $F:E\times\Gamma\to\mathbb{C}$ and their evolution
$F\mapsto F(Z_t)=F\circ\alpha^t(Z_0)$.

In order to show that every state space $E$ admits a universal dilation, now we consider a particular
class of universal dilations. Let us describe all the special requirements we are interested in. 

First of all we want the sample space $\Omega$ to be just $E\times\Gamma$, the state space of the 
global system. As we want it to describe all the possible initial global states, we ask
the random variable $Z_0$ to be the identity function and $X_0$ and $Y_0$ to be the coordinate 
variables: if $\omega=(i,\gamma)$, then $Z_0(\omega)=\omega$, $X_0(\omega)=i$ and
$\Upsilon_0(\omega)=\gamma$. Thus, for all $t\geq1$, $Z_t=(X_t,\Upsilon_t)=\alpha^t\circ
Z_0=Z_0\circ\alpha^t$, $X_t=X_0\circ\alpha^t$ and $\Upsilon_t=\Upsilon_0\circ\alpha^t$. 

We are interested in an
environment $\Gamma=G^\mathbb{Z}=\bigtimes_{n\in\mathbb{Z}} G_n$ with $\quad G_n=G$ 
finite set. In this case the environment state $\gamma\in\Gamma$ has infinite components $g_n\in G_n$,
$n\in\mathbb{Z}$, and we introduce also the coordinate variables $Y_m(\omega)=g_m$, the random $m$-th
components of the environment. Then, considered the power $\sigma$-field $\mathcal{E}$ on $E$ and 
the power $\sigma$-field on $G$, we ask $\mathcal{G}$ to be the cylindric $\sigma$-field on $\Gamma$, 
and $\mathcal{F}$ to be 
$\mathcal{E}\otimes\mathcal{G}=\sigma(X_0,\Upsilon_0)=\sigma(X_0,Y_n,\,n\in\mathbb{Z})$
on $\Omega$, so that $X_0$, $\Upsilon_0$ and the $Y_n$ are all measurable. 

Supposing that at time 0
only $X_0$ is observed and that at each following instant $t$ only the information carried by $Y_t$ is
acquired, we want the filtration $\mathcal{F}_0=\sigma(X_0)$, 
$\mathcal{F}_t=\sigma(X_0, Y_s ,1\leq s\leq t)$, for $t\geq1$. Coherently we also want an evolution
$\alpha$ such that each $Y_t$ is involved in the interaction with the system only once, between time
$t-1$ and time $t$, so that actually $X_t=b(X_{t-1},Y_t)$ for a suitable function $b:E\times G\to E$.
In this case $X_t$ is automatically adapted to $\mathcal{F}_t$.
This is possible if we introduce an invertible map $\varphi_1$ on $E\times G_1$, denoted by
$\varphi_1(i,g)=\big(\varphi_1^E(i,g),\varphi_1^G(i,g)\big)$, and the invertible shift $\vartheta$ on
$\Gamma=\bigtimes_{n\in\mathbb{Z}} G_n$ mapping each $\gamma=(g_n)_n$ to the environment state
$\vartheta(\gamma)$ with each $n$-th component equal to $g_{n+1}$. Then, trivially extended 
$\varphi_1$ and
$\vartheta$ to $\Omega$, it is enough to set $\alpha=\vartheta\circ\varphi_1$, getting
$b=\varphi_1^E$. Indeed, when $\alpha$ is applied for the first time between time 0 and time 1, 
the map $\varphi_1$ couples the system state $X_0$ with $Y_1$, giving the new system state
$X_1=X_0\circ\alpha=X_0\circ\varphi_1=\varphi_1^E(X_0,Y_1)$, and then
the shift $\vartheta$ prepares $Y_2$ for the following interaction 
with $X_1$. Thus $\vartheta$ could be interpreted as a
free evolution of the environment. 
If we explicitly introduce also the random variables $Y^{(t)}_n=Y_n\circ\alpha^t$, the
$n$-th environment components at time $t$, then at time 1 the environment state is 
$\Upsilon_1=(Y^{(1)}_n)_n$ with $Y^{(1)}_0=\varphi^G_1(X_0,Y_1)$ and with
$Y^{(1)}_n=Y_{n+1}$ for all $n\neq0$. 
After $t$ steps, $X_t=\varphi_1^E(X_{t-1},Y_t)=\varphi_1^E(X_{t-1},Y^{(t-1)}_1)$,
while $Y^{(t)}_0=\varphi^G_1(X_{t-1},Y^{(t-1)}_1)$ and $Y^{(t)}_n=Y^{(t-1)}_{n+1}$ for all $n\neq0$.

Finally, we require the probabilities $\mathbb{P}_{k,\{P\}}$ to be $\delta_k\otimes Q_{\{P\}}$ on
$\mathcal{F}=\mathcal{E}\otimes\mathcal{G}$, with $Q_{\{P\}}$, the distribution of $\Upsilon_0$, 
such that $(Y_n)_{n\leq0}$, $Y_1$, $Y_2, \ldots$ are independent. Then this independence
guarantees the Markov property with respect to $\mathcal{F}_t$ for the process $X_t$.
 
A universal dilation like this will be called \textbf{standard} in the following.
Summarizing, a universal dilation is standard if 
\begin{itemize}
\item $\Gamma=G^\mathbb{Z}=\bigtimes_{n\in\mathbb{Z}} G_n$, $\quad G_n=G$ finite set with power
$\sigma$-field, $\quad \mathcal{G}=$ cylindric $\sigma$-field on $\Gamma$,
\item $\Omega=E\times\Gamma, \qquad \omega=(i,\gamma)=(i,(g_n)_n)\in\Omega, 
\quad i\in E, \quad \gamma\in\Gamma, \quad g_n\in G$,
\item $X_0(\omega)=i, \quad\; \Upsilon_0=(Y_n)_{n\in\mathbb{Z}}, \quad\; Y_n:\Omega\to G,
\quad\; \Upsilon_0(\omega)=\gamma, \quad\; Y_m(\omega)=g_m, \quad\; Z_0(\omega)=\omega$,
\item $\mathcal{F}=\sigma(X_0,\Upsilon_0), \qquad \mathcal{F}_t=\sigma(X_0, Y_s ,1\leq s\leq t),
\quad t\geq0$,
\item $\mathbb{P}_{k,\{P\}}=\delta_k\otimes Q_{\{P\}}$, with $Q_{\{P\}}$ such that $(Y_n)_{n\leq0}$,
$Y_1$, $Y_2, \ldots$ are independent,
\item $\alpha=\vartheta\circ\varphi_1$, with an invertible map $\varphi_1:E\times G_1\to E\times G_1$,
where $\varphi_1(i,g)=\big(\varphi_1^E(i,g),\varphi_1^G(i,g)\big)$,
and the invertible shift $\vartheta:\Gamma\to \Gamma$,
\begin{equation*}
    (g_n)_{n\in\mathbb{Z}}\in\Gamma,\quad g_n\in G_n \quad \mapsto 
    \quad (g_{n+1})_{n\in\mathbb{Z}}\in\Gamma,\quad g_{n+1}\in G_n.
\end{equation*}
\end{itemize} 
Then, we have
\begin{itemize}
\item $Z_t=\alpha^t\circ Z_0=Z_0\circ\alpha^t=(X_t,\Upsilon_t)$, 
\item $\Upsilon_t=\Upsilon_0\circ\alpha^t=(Y_n^{(t)})_{n\in\mathbb{Z}}$, 

$\displaystyle Y^{(t)}_n = Y_n\circ\alpha^t = \begin{cases} 
    Y_{n+1}^{(t-1)},& n\neq0,\\ 
    \varphi^G_1(X_{t-1},Y_1^{(t-1)}),& n=0,
\end{cases} = \begin{cases} 
    Y_{n+t},& n\leq-t,\;n\geq1\\ 
    \varphi^G_1(X_{t-1+n},Y_{n+t}),& -t+1\leq n\leq0.
\end{cases}$
\item $X_t=X_0\circ\alpha^t=\varphi_1^E(X_{t-1},Y_t)=\varphi_1^E(X_{t-1},Y^{(t-1)}_1)$, 
\end{itemize}

A standard universal dilation is therefore specified by the term
$\big(G^\mathbb{Z},\mathcal{G},\alpha,(Q_{\{P\}})_{\{P\}\in\mathcal{P}}\big)$, where
$(G^\mathbb{Z},\mathcal{G})$ represents the environment
state space, $\alpha$ the one-step evolution of the system with the environment, and $Q_{\{P\}}$ is a
family of
initial distributions for the environment state which give rise to the corresponding Markov
evolutions for the system. 

With a standard universal dilation the evolution of every global random variable
$F\in\mathcal{L}^\infty(\mathcal{F})$ is described by the group of $*$-automorphisms $F\mapsto
F\circ\alpha^t$. 

Note that a standard universal dilation always allows to see every Markov chain also as an
automaton with independent random inputs $Y_t$, which now are provided by the environment via the
dynamics $\alpha$.

A standard universal dilation is not uniquely determined by the state space
$E$, but it always exists.

\begin{theorem}\label{cud}
For every finite state space $E$, there exists a standard universal dilation\\
$\big(G^\mathbb{Z},\mathcal{G},\alpha,(Q_{\{P\}})_{\{P\}\in\mathcal{P}}\big)$ of the Markov 
evolutions on $E$.
\end{theorem}

\noindent {\sl Proof.}
We have to exhibit a proper set $G$, together with the coupling $\varphi_1$ and the probability 
measures $\mathbb{P}_{k,\{P\}}$.

Given $E=\{1,\ldots,N\}$ and the set $L$ labelling the all possible maps $\beta:E\to E$, we set
\begin{equation*}
    G=E\times L, \qquad\qquad
    i,j,k\in E, \qquad \ell\in L, \qquad g=(j,\ell)\in G.
\end{equation*}

Now, taken two points
$\big(i,(1,\ell)\big)\neq\big(i',(1,\ell')\big)$ in $E\times G$, we get
$\big(\beta_\ell(i),(i,\ell)\big)\neq\big(\beta_{\ell'}(i'),(i',\ell')\big)$ 
and so we can find an invertible map
\begin{equation}\label{coupling}
    \varphi:E\times G\to E\times G, \qquad\qquad 
    \varphi\big(i,(j,\ell)\big) = \begin{cases}\big(\beta_\ell(i),(i,\ell)\big), \quad & \text{if
    }j=1,\\ \quad \ldots \;, & \text{if }j\neq1.\end{cases}
\end{equation}
The coupling $\varphi_1$ is defined to be $\varphi$ when $G=G_1$.

Given $\{P\}$, 
we choose a decomposition \eqref{dcc} for every $P(t)$, thus obtaining the distributions $q(t)$ on 
$L$;
fixed an arbitrary state in $E$, say 1, for every distribution $q$ on $L$ we define the distribution
$\widetilde q=\delta_1\otimes q$ on $G$, that is $\widetilde q_g=\widetilde q_{(j,\ell)}=
\delta_{j,1}\, q_\ell$. Fixed also an arbitrary distribution $Q_0$ for $\bigtimes_{n\leq0} G_n$,
we define $Q_{\{P\}}=Q_0 \otimes \big(\bigotimes_{t\in\mathbb{N}} \widetilde q(t)\big)$ on
$\mathcal{G}$ and, given $k\in E$, we define  
$\mathbb{P}_{k,\{P\}}=\delta_k\otimes Q_{\{P\}}$ on $\mathcal{F}$, that is
\begin{equation*}
    \mathbb{P}_{k,\{P\}}\Big(X_0=i,Y_{-m}=g_{-m},\ldots,Y_{m}=g_m\Big)=
    \delta_{i,k} \, Q_0(g_{-m},\ldots,g_0) \, \prod_{n=1}^m\widetilde q_{g_n}(n), \qquad \forall
    m\geq1.
\end{equation*}
Then $X_0=k$ a.s. and $\Upsilon_0$ has distribution $Q_{\{P\}}$, so that $(Y_n)_{n\leq0}\sim Q_0$ 
while $Y_n\sim\widetilde q(n)$ for every $n\geq1$ and they are all independent.
Of course, if we change $k$ and $\{P\}$, then $X_0\sim\delta_k$ for every $\{P\}$ and $\Upsilon_0\sim
Q_{\{P\}}$ for every $k$. Moreover, every 
$\big(\Omega, \mathcal{F}, (\mathcal{F}_t)_{t\geq0}, (X_t)_{t\geq0}, (\mathbb{P}_{k,\{P\}})_{k\in
E}\big)$ is a Markov chain with transition matrices $\{P\}$, independently of $Q_0$ and 
of the definition of $\varphi\big(i,(j,\ell)\big)$ for $j\neq1$. Indeed, for every $t\geq1$ and 
every event $(Y_t=g_t,\ldots,Y_1=g_1,X_0=i_0)$ with positive probability $\mathbb{P}_{k,\{P\}}$, 
that implies $i_0=k$ and $g_s=(1,\ell_s)$ for all $1\leq s\leq t$, it holds
\begin{equation*}\begin{split}
    \mathbb{P}&_{k,\{P\}}\Big(X_{t+1}=j\Big|Y_t=(1,\ell_t),\ldots,Y_1=(1,\ell_1),X_0=k\Big)\\
    &=\sum_{\ell\in L}\mathbb{P}_{k,\{P\}}\Big(X_{t+1}=j \Big| 
    Y_{t+1}=(1,\ell),Y_t=(1,\ell_t),\ldots,Y_1=(1,\ell_1),X_0=k\Big) \\
    &{}\qquad\qquad\qquad\qquad\qquad\qquad\qquad
    {}\times\mathbb{P}_{k,\{P\}}\Big(Y_{t+1}=(1,\ell)\Big|Y_t=(1,\ell_t),\ldots,Y_1=(1,\ell_1),X_0=k\Big)\\
    &=\sum_{\ell\in L}\delta_{\beta_\ell\circ\beta_{\ell_t}\circ\cdots\circ\beta_{\ell_1}(k),j}\;
    \;\mathbb{P}_{k,\{P\}}\Big(Y_{t+1}=(1,\ell)\Big)\\
    &=\sum_{\ell\in L}\delta_{\beta_\ell(\beta_{\ell_t}\circ\cdots\circ\beta_{\ell_1}(k)),j} 
    \;\;q_\ell(t+1) \\
    &= p_{\beta_{\ell_t}\circ\cdots\circ\beta_{\ell_1}(k)\,j}(t+1),
\end{split}\end{equation*}
where the last equality follows from \eqref{detmatrix} and \eqref{dcc}. Thus
\begin{equation*}
    \mathbb{P}_{k,\{P\}}(X_{t+1}=j|\mathcal{F}_t) = p_{X_tj}(t+1) = \mathbb{P}_k(X_{t+1}=j|X_t), 
    \qquad \forall \{P\},k,j,t.
\end{equation*}
Therefore 
$\big(G^\mathbb{Z},\mathcal{G},\alpha,(Q_{\{P\}})_{\{P\}\in\mathcal{P}}\big)$
is a standard universal dilation
of the Markov evolutions on $E$. 
\QED

In this construction, even if each $g\in G$
has two components, $g=(j,\ell)$, the probability is always
concentrated only  on those $g$ of the kind
$g=(1,\ell)$, but we need the first component $j$ to define an
invertible $\varphi$. Analogously, we are considering the evolution only for positive
times so that all the components $g_n$, $n\leq0$, are never involved in the interaction with the
system, but they are needed to define an invertible shift $\vartheta$.

\bigskip

\noindent \textbf{Example 1.} As a first example, let us characterize the dilation provided by Theorem
\ref{cud} in the degenerate case of a homogeneous deterministic evolution. Let every $P(t)$ be
equal to the deterministic matrix $D_{\bar\ell}$ associated to the map $\beta_{\bar\ell}$ on $E$. 
Then, for positive $t$, the distributions $q(t)$ must equal $\delta_{\bar\ell}$, so that
$\widetilde{q}(t)=\delta_1\otimes\delta_{\bar\ell}=\delta_{(1,{\bar\ell})}$ and the environment
components $Y_t$, the ones which will interact with the system, turn out to be deterministic. For the
sake of simplicity, let us choose the same degenerate distributions also for all $t\leq0$, that is
$Q_0=\bigotimes_{t\leq0}\delta_{(1,{\bar\ell})}$, and, for every system initial state $k$, let
$\mathbb{P}_{k,D_{\bar\ell}}=\delta_k\otimes(\bigotimes_\mathbb{Z}\delta_{(1,{\bar\ell})})$ denote 
the corresponding probability on $\mathcal{F}$. Then, under every $\mathbb{P}_{k,D_{\bar\ell}}$, the 
environment initial state $\Upsilon_0$ itself is deterministic and, for every coupling $\varphi$ 
satisfying \eqref{coupling}, the following equalities hold almost surely
\begin{gather*}
    X_0=k, \qquad \Upsilon_0=((1,{\bar\ell})_n),\qquad Z_0=(k,(1,{\bar\ell})_n), \\
    Z_t=\alpha^t \circ Z_0=({\beta_{\bar\ell}}^t(k),(g_n)_n), \;\; \forall t\geq0, \qquad
    \text{where } g_n = \begin{cases}
    (1,{\bar\ell}),&n\leq-t,\;n\geq1,\\({\beta_{\bar\ell}}^{t-1+n}(k),{\bar\ell}),&-t+1\leq n\leq0,
    \end{cases}\\
    X_t={\beta_{\bar\ell}}^t(k).
\end{gather*}
Note that, even if $D_{\bar\ell}$ is a permutation, that is even if $\beta_{\bar\ell}$ is invertible
(and no dilation would be needed), the dilation provided by Theorem \ref{cud} is not trivial at all.
Indeed, just because of the universality of the construction, a dilation
$\big(G^\mathbb{Z},\mathcal{G},\alpha,Q_{\{P\}}\big)$ is usually non minimal for a particular 
evolution $\{P\}$.

\bigskip

\noindent \textbf{Example 2.} As a second example, let us consider a homogeneous evolution in
$E=\{1,2\}$, say $P(t)=P$ for every $t$. In this case there are four deterministic matrices
\begin{equation*}
    D_1=\begin{pmatrix}1&0\\1&0\end{pmatrix}, \qquad
    D_2=\begin{pmatrix}1&0\\0&1\end{pmatrix}, \qquad
    D_3=\begin{pmatrix}0&1\\1&0\end{pmatrix}, \qquad
    D_4=\begin{pmatrix}0&1\\0&1\end{pmatrix},
\end{equation*}
and we decompose the stochastic matrix $P$ as
\begin{equation*}
    P=\begin{pmatrix}p_{11}&p_{12}\\p_{21}&p_{22}\end{pmatrix}=
    p_{11}p_{21}\,D_1 + p_{11}p_{22}\,D_2 + p_{12}p_{21}\,D_3 + p_{12}p_{22}\,D_4 = 
    q_1\,D_1 + q_2\,D_2 + q_3\,D_3 + q_4\,D_4.
\end{equation*}
Therefore $L=\{1,2,3,4\}$ and the environment state space is $\Gamma=G^\mathbb{Z}$ with
$G=\{1,2\}\times\{1,2,3,4\}$. The distribution $q$ on $L$ gives the distribution
$\widetilde{q}=\delta_1\otimes q$
on $G$ and the distribution $\bigotimes_\mathbb{Z}\widetilde{q}$ on $G^\mathbb{Z}$. Choosing this
latter as the distribution $Q_P$, under every
$\mathbb{P}_{k,P}=\delta_k\otimes(\bigotimes_\mathbb{Z}\widetilde{q})$ 
the initial states of the system and of the environment are
\begin{gather*}
    X_0=k \;\; \text{a.s.}, \qquad\qquad \Upsilon_0=(Y_n)_n, \quad Y_n
    \;\;\text{i.i.d.}, \quad Y_n\sim \widetilde{q},
\end{gather*}
so that, for example, $\mathbb{P}_{k,P}(Y_n=(1,3))=p_{12}p_{21}$.
Even if the definition of the coupling \eqref{coupling} for $j\neq1$ plays no role, 
let us fix $\varphi$ setting
\begin{equation*}
    \varphi:E\times G\to E\times G, \qquad\qquad 
    \varphi\big(i,(j,\ell)\big) := \begin{cases}
    \big(\beta_\ell(i),(i,\ell)\big), \quad & \text{if }j=1,\\
    \big(\beta_{\ell+2}(i),(i,\ell)\big),\quad & \text{if }j=2,\,i=1,\\ 
    \big(\beta_{\ell+1}(i),(i,\ell)\big),\quad & \text{if }j=2,\,i=2,\end{cases}
\end{equation*}
where $\ell+1$ and $\ell+2$ are sums modulo 4. 
Then, always under $\mathbb{P}_{k,P}$, after one step 
the states of the systems are both random and, precisely,
\begin{equation*}
    X_1=\varphi_1^E(k,Y_1), \qquad\qquad \Upsilon_1=(Y^{(1)}_n)_n, \quad 
    Y^{(1)}_n=\begin{cases}
    Y_{n+1}, \quad & n\leq-1,\;n\geq1\\ \varphi^G_1(k,Y_1),\quad & n=0,\end{cases}
\end{equation*}
where $X_1=\varphi_1^E(k,Y_1)=\beta_\ell(k)$ when $Y_1=(1,\ell)$. Therefore
\begin{multline*}
    \mathbb{P}_{k,P}(X_1=j)
    =\sum_{\ell=1}^4\mathbb{P}_{k,P}(X_1=j|Y_1=(1,\ell)) \, \mathbb{P}_{k,P}(Y_1=(1,\ell))
    =\sum_{\ell=1}^4 \delta_{\beta_\ell(k),j} \, q_\ell
    =\sum_{\ell=1}^4 (D_\ell)_{k,j} \, q_\ell \\ = p_{kj}.
\end{multline*}
Similar calculations for $t\geq2$ then prove that 
$\big(\Omega, \mathcal{F}, (\mathcal{F}_t)_{t\geq0},(X_t)_{t\geq0}, (\mathbb{P}_{k,P})_{k\in E}\big)$
is a homogeneous Markov chain with transition matrix $P$.

\bigskip

Let us remark that, as long as we consider only system random variables neglecting the 
environment, we can avoid the shift $\vartheta$ 
and, starting from this dilation, we 
can define a deterministic, invertible, but inhomogeneous global evolution 
$\widetilde\varphi_t$ which never involves the environment
components $g_n$ for $n\leq0$, but which generates the
same Markov evolutions for the system. Consider indeed
the invertible maps on $\Omega$
\begin{equation}\label{cocycle}
    \varphi_t:=\vartheta^{-(t-1)} \circ \varphi_1 \circ \vartheta^{t-1}, \qquad
    \widetilde\varphi_t:=\varphi_t\circ\cdots\circ\varphi_1=\vartheta^{-t}\circ\alpha^t, \qquad
    t\geq1,
\end{equation}
where $\varphi_t$ actually acts only on $E\times G_t$, where it is just the coupling \eqref{coupling},
while $\widetilde\varphi_t$ actually acts only on $E\times G_1\times\cdots\times G_t$. Then, for every
$t\geq0$, we get $X_t=X_0\circ\alpha^t=X_0\circ\widetilde\varphi_t$, so that the same system 
evolution is
obtained without involving $g_n$ for $n\leq0$. In particular the stochastic evolution of a system
random variable $f\in\mathcal{L}^\infty(E)$ is given by the $*$-unital homomorphism
\begin{equation}\label{clfl}
    j_t:\mathcal{L}^\infty(E)\to\mathcal{L}^\infty(\mathcal{F}_t), \qquad f\mapsto
    j_t[f]:=f(X_t)=f\circ\widetilde\varphi_t, \qquad t\geq0,
\end{equation}
which is injective as $\widetilde\varphi_t$ is invertible.
If $\mathbb{E}_g[f\circ\varphi]\in\mathcal{L}^\infty(E)$ is the system random variable $i\mapsto
f\circ\varphi(i,g)$, then the stochastic evolution \eqref{clfl} satisfies
\begin{equation}\label{clfleq}
    j_0[f]=f(X_0), \qquad j_t[f]=\sum_{g\in G}
    j_{t-1}\Big[\mathbb{E}_g[f\circ\varphi]\Big]\,I_{(Y_t=g)}, \qquad
    \forall f\in\mathcal{L}^\infty(E), \; t\geq1,
\end{equation}
where $j_{t-1}\Big[\mathbb{E}_g[f\circ\varphi]\Big]$ are $\mathcal{F}_{t-1}$-measurable
and $I_{(Y_t=g)}$ are the indicators of the events $(Y_t=g)$.
And now we could even reduce the sample space $\Omega$
from $E\times G^\mathbb{Z}$ to $E\times G^\mathbb{N}$, restricting here $\mathcal{F}$,
$\mathcal{F}_t$ and
$\mathbb{P}_{k,\{P\}}$. This alternative construction has its counterpart in Quantum Probability and,
analogously to what happens in the Quantum framework, one could start with it and then recover the
whole
dilation $\big(G^\mathbb{Z},\mathcal{G},\alpha,(Q_{\{P\}})_{\{P\}\in\mathcal{P}}\big)$,
so that the two constructions can be considered equivalent descriptions of the same situation.

To conclude this section, let us add two remarks about the choice of probabilities $\mathbb{P}$ 
different from the selected $\mathbb{P}_{k,\{P\}}$, again in the case of 
the dilation built in Theorem \ref{cud}.
Firstly, the decomposition \eqref{dcc} is not unique, just as the choice of
$Q_0$, so that we could find other probabilities $\mathbb{P}$ on $\mathcal{F}$ inducing the same 
Markov evolutions for the system. Secondly, we could also consider non-factorized $Q$, or even 
non-factorized $\mathbb{P}$, thus obtaining new stochastic
processes $\big(\Omega, \mathcal{F}, (\mathcal{F}_t)_{t\geq0}, (X_t)_{t\geq0}, \mathbb{P}\big)$, 
without any Markov property guaranteed, which would depend also on the definition of $\varphi$ for
$j\neq1$.

\section{UNIVERSAL DILATION OF MARKOV SEMIGROUPS AND DILATIONS OF QUANTUM DYNAMICAL SEMIGROUPS}

We want to compare the dilation built in  Theorem \ref{cud} with similar constructions typical in
Quantum Probability: given a stochastic matrix $P$ on $E$, we show that the classical evolution 
$(P^t)_{t\geq0}$ (Classical Markov Semigroup, CMS) can be 
extended to a quantum evolution $(T^t)_{t\geq0}$ (Quantum Dynamical Semigroup, QDS) 
for which we can find a quantum dilation which is itself an extension of the classical
dilation $\big(G^\mathbb{Z},\mathcal{G},\alpha,Q_{P}\big)$. Of course, $Q_P$ and $\mathbb{P}_{k,P}$ 
stay for $Q_{\{P\}}$ and $\mathbb{P}_{k,\{P\}}$ when the sequence $\{P\}$ is constant.

Given a complex separable Hilbert space $\is$, let us denote its vectors by $h$, or $|h\rangle$ using
Dirac's notation, so that $\langle h'|h\rangle$ denotes the scalar product (linear in $h$) and
$|h'\rangle\langle h|$ denotes the bounded operator $h'' \mapsto \langle h|h''\rangle \,h'$. Let
$\mathcal{B}(\is)$ be the complex $*$-algebra of bounded operators on $\is$ and let $\mathcal{S}(\is)$ be the
convex set of states on $\mathcal{B}(\is)$,
$\mathcal{S}(\is)=\Big\{\rho\in\mathcal{B}(\is) \Big| \rho\geq0, \tr\rho=1\Big\}$, with
its extreme points (pure states) given by the rank-1 projections
$|h\rangle\langle h|$, $\|h\|=1$, so that for every $\rho$ in
$\mathcal{S}(\is)$ there exists an orthonormal basis $\{h_i\}$ such that $\rho = \sum_i
p_i\,|h_i\rangle\langle h_i|$, $p_i\geq0$, $\sum_i p_i=1$.
For every $\rho$ in $\mathcal{S}(\is)$, the term $\big(\mathcal{B}(\is),\rho\big)$ is a quantum
probability space and $\tr[a\rho]$ is the expected value of $a\in\mathcal{B}(\is)$.

A family $\{a_\lambda\}_{\lambda\in\Lambda}\subseteq\mathcal{B}(\is)$ of normal 
operators is commuting if and only if the operators $a_\lambda$ admit a common spectral
representation, i.e.\ a classical
probability space $(\Omega,\mathcal{F},\mathbb{P})$, a unitary operator $u:\is\to L^2(\mathbb{P})$ and
a family $\{f_\lambda\}_{\lambda\in\Lambda}\subseteq L^\infty(\mathbb{P})$ of random variables on
$\Omega$, such that
\begin{equation*}
    a_\lambda = u^{-1}\,m_{f_\lambda}\,u, \qquad \forall\lambda\in\Lambda,
\end{equation*}
where, for every $f$ in $L^\infty(\mathbb{P})$, $m_f$ denotes the bounded multiplication operator by
$f$ on $L^2(\mathbb{P})$. Thanks to $u$ the operators $a_\lambda$ admit a common probabilistic
interpretation as the random variables $f_\lambda$, with a
joint distribution which depends on the state $\rho\in\mathcal{S}(\is)$: 
via the spectral representation, every $\rho= \sum_i p_i\,|h_i\rangle\langle h_i|$ 
defines a probability $\mathbb{P}_\rho$ on $(\Omega,\mathcal{F})$,
\begin{equation*}
    \frac{\dd\mathbb{P}_\rho}{\dd\mathbb{P}} = \sum_i p_i\,|u\,h_i|^2,
\end{equation*}
and so a joint distribution for $\{f_\lambda\}_{\lambda\in\Lambda}$ and hence for
$\{a_\lambda\}_{\lambda\in\Lambda}$,
\begin{equation*}
    \mathbb{E}_\rho\big[\eta(f_{\lambda_1},\ldots,f_{\lambda_n})\big] = 
    \int_\Omega \eta(f_{\lambda_1},\ldots,f_{\lambda_n})\,\dd \mathbb{P}_\rho = 
    \tr\big[\eta(a_{\lambda_1},\ldots,a_{\lambda_n})\,\rho\big],
\end{equation*}
$n\in\mathbb{N}$, $\{\lambda_1,\ldots,\lambda_n\}\subseteq\Lambda$, $\eta:\mathbb{C}^n\to\mathbb{R}$
bounded and continuous.

An arbitrary one-step evolution in $\mathcal{B}(\is)$ is given by a stochastic map $T$, that is 
a completely positive, identity preserving, linear 
operator $T:\mathcal{B}(\is)\to\mathcal{B}(\is)$. A stochastic map $T$ is 
invertible, i.e.\ there exists another stochastic map $T'$ such that $T\circ T'=T'\circ
T=\operatorname{Id}_{\mathcal{B}(\is)}$, if and only if $T[a]=u^*au$ with $u$ unitary operator
on $\is$ ($T$ $*$-automorphism of $\mathcal{B}(\is)$). 
A homogeneous evolution in discrete time is given by a QDS, that is
a semigroup $(T^t)_{t\geq0}$ with $T$ as above.

Let us now give first the definition of Quantum Stochastic Flow (QSF), then the definitions
of dilation of a stochastic map and of a QDS, and finally the definitions of extension of a CMS and of
a classical standard dilation.

Given two Hilbert spaces $\is$ and $\mathcal{K}$ and a trace class operator 
$\sigma$ on $\mathcal{K}$, the
\textbf{conditional expectation} of an operator $A$ in $\mathcal{B}(\is\otimes\mathcal{K})=
\mathcal{B}(\is)\otimes\mathcal{B}(\mathcal{K})$ with respect to $\sigma$ is the operator
$\mathbb{E}_\sigma[A]$ in $\mathcal{B}(\is)$ defined by
\begin{equation*}
    \tr\Big[\mathbb{E}_\sigma[A]\cdot\rho\Big]=\tr\big[A\cdot\rho\otimes\sigma\big],
    \qquad \forall\rho\in\mathcal{S}(\is).
\end{equation*}

Given $\is$ and infinitely many copies $\mathfrak{Z}_n$ of a Hilbert space $\mathfrak{Z}$, introduced
the increasing sequence of algebras
\begin{equation*}
    \mathcal{B}_{(0,t]}=\mathcal{B}(\is)\otimes\mathcal{B}\Big(\bigotimes_{n=1}^t\mathfrak{Z}_n\Big),
    \qquad t\geq1,
\end{equation*}
a \textbf{QSF} on $\mathcal{B}(\is)$ in discrete time is a family of $*$-unital 
homomorphisms $j_t:\mathcal{B}(\is)\to\mathcal{B}_{(0,t]}$, $t\geq0$, induced by a $*$-unital 
homomorphism $\nu:\mathcal{B}(\is)\to\mathcal{B}(\is\otimes\mathfrak{Z})$ through the equations
\begin{equation}\label{qsf}
    j_0=\operatorname{Id}_{\mathcal{B}(\is)}, \qquad
    j_t[a]=\sum_{zz'}j_{t-1}\Big[\mathbb{E}_{|z\rangle\langle z'|}\big[\nu(a)\big]\Big] \otimes
    |z'\rangle\langle z|, \quad t\geq1,
\end{equation}
where $\{|z\rangle\}$ is a given basis of $\mathfrak{Z}$.

Given a Hilbert space
$\mathfrak{Z}$, a unitary operator $V$ on $\is\otimes\mathfrak{Z}$ and a state 
$\sigma\in\mathcal{S}(\mathfrak{Z})$, the term $(\mathfrak{Z},V,\sigma)$ is a \textbf{dilation of a 
stochastic map} $T$ on $\mathcal{B}(\is)$ if
\begin{equation*}
    T[a]=\mathbb{E}_\sigma\big[V^*\cdot a\otimes\bbbone_\mathfrak{Z}\cdot V\big], \qquad
    \forall a\in\mathcal{B}(\is),
\end{equation*}
so that the one-step evolution $T$ of the quantum system $\is$ is represented in
terms of an invertible evolution $V$ of the system coupled with an environment $\mathfrak{Z}$ 
with initial state $\sigma$. A dilation always exists, not unique, and it can always be chosen with 
a pure state $\sigma=|\phi\rangle\langle\phi|$, $\phi\in\mathfrak{Z}$.

Given a Hilbert space $\mathcal{K}$, a unitary operator $U$ on $\is\otimes\mathcal{K}$ and a state 
$\sigma\in\mathcal{S}(\mathcal{K})$, the term $(\mathcal{K},U,\sigma)$ is a \textbf{dilation of a QDS} 
$T^t$ on $\mathcal{B}(\is)$ by the group of $*$-automorphisms $A\mapsto {U^*}^tAU^t$ on
$\mathcal{B}(\is)\otimes\mathcal{B}(\mathcal{K})$ if
\begin{equation*}
    T^t[a]=\mathbb{E}_\sigma\big[{U^*}^t\cdot a\otimes\bbbone_\mathfrak{Z}\cdot U^t\big], \qquad
    \forall a\in\mathcal{B}(\is), \;t\geq0,
\end{equation*}
so that the discrete time evolution $T^t$ of the quantum system $\is$ is represented in
terms of an invertible homogeneous evolution $U^t$ of the system coupled with an environment 
$\mathcal{K}$ with initial state $\sigma$.

Given a Hilbert space $\mathfrak{Z}$, a unitary operator $V$ on 
$\is\otimes\mathfrak{Z}$ and a unit vector $\phi$ in $\mathfrak{Z}$, the
term $(\mathfrak{Z},V,\phi)$ gives a \textbf{dilation of a QDS} $T^t$ on $\mathcal{B}(\is)$ 
by the QSF
$(j_t)_{t\geq0}$ induced by the $*$-unital homomorphism 
$\theta:\mathcal{B}(\is)\to\mathcal{B}(\is\otimes\mathfrak{Z})$, 
$\nu[a]=V^*\cdot a\otimes\bbbone_\mathfrak{Z}\cdot V$ if
\begin{equation*}
    T^t[a]=\mathbb{E}_{\phi^{\otimes t}}\big[j_t[a]\big], \qquad
    \forall a\in\mathcal{B}(\is), \;t\geq0,
\end{equation*}
so that the discrete time evolution $T^t$ of the quantum system $\is$ is represented as the
expectation of a QSF associated to an invertible coupling $V$, 
analogously to the stochastic evolution \eqref{clfl}, \eqref{clfleq}.

A typical construction employs the same dilation $(\mathfrak{Z},V,\phi)$ of a stochastic map $T$ 
to dilate the QDS $T^t$ by the associated QSF and, at the same time, by a related group of 
$*$-automorphisms. Taken infinitely many copies $\mathfrak{Z}_n$ of $\mathfrak{Z}$,
$n\in\mathbb{Z}$, let us define the Hilbert space
$\mathcal{K}=\bigotimes_{n\in\mathbb{Z}}\mathfrak{Z}_n$, infinite tensor product 
(\cite{P} Ex.~15.10 or \cite{KR} Ex.~11.5.29)
with respect to a stabilizing sequence of norm-1 vectors
$\psi_n\in\mathfrak{Z}_n$ such that $\psi_n=\phi$ for every $n\geq1$. 
Let us set $\Psi=\bigotimes_{n\in\mathbb{Z}}\psi_n$, norm-1 vector in $\mathcal{K}$.
Let us also embed each $\mathcal{B}_{(0,t]}$ into
$\mathcal{B}(\is)\otimes\mathcal{B}(\mathcal{K})$ by tensorizing with the identity. In particular let
us denote by $V_n$ the unitary operator $V$ on $\is\otimes\mathfrak{Z}_n$, let us identify all of them
with their extension to $\is\otimes\mathcal{K}$, and, for every $t\geq1$, let us define
$\widetilde{V}_t=V_t \cdots V_1$, unitary operator acting on
$\is\otimes\mathfrak{Z}_1\otimes\cdots\otimes\mathfrak{Z}_t$. Then $j_t[a]:=\widetilde{V}_t^*\cdot
a\otimes\bbbone_\mathcal{K} \cdot\widetilde{V}_t$ is the QSF satisfying
\begin{equation*}
    j_t[a]=\sum_{zz'}j_{t-1}\Big[\mathbb{E}_{|z\rangle\langle z'|}
    \big[V^*\cdot a\otimes\bbbone_\mathfrak{Z}\cdot V\big]\Big] \otimes
    |z'\rangle\langle z|, \qquad t\geq1,
\end{equation*}
for every basis $\{|z\rangle\}$ in $\mathfrak{Z}$, and it dilates the QDS $T^t$:
\begin{equation*}
    T^t[a]=\mathbb{E}_\Psi\big[j_t[a]\big], \qquad\forall a\in\mathcal{B}(\is),\;t\geq0.
\end{equation*}
Moreover, introduced the shift operator $\Theta:\mathcal{K}\to\mathcal{K}$
\begin{equation}\label{qshift}
    \textstyle \bigotimes_n z_n\in\mathcal{K},\quad z_n\in\mathfrak{Z}_n \quad \mapsto 
    \quad \bigotimes_n z_{n+1}\in\mathcal{K},\quad z_{n+1}\in\mathfrak{Z}_n
\end{equation}
(where, of course, $z_n=\phi$ for large $n$), extended it to a unitary operator
$\Theta:\is\otimes\mathcal{K}\to\is\otimes\mathcal{K}$, let us define the unitary operator
\begin{equation}\label{glev}
    U= \Theta\,V_1.
\end{equation}
Then $(\mathcal{K},U,\Psi)$ defines a group of $*$-automorphisms of
$\mathcal{B}(\is)\otimes\mathcal{B}(\mathcal{K})$ 
which dilates the QDS $T^t$:
\begin{equation*}
    T^t[a]=\mathbb{E}_\Psi\big[{U^*}^t\cdot a\otimes\bbbone_\mathcal{K}\cdot U^t\big], 
    \qquad\forall a\in\mathcal{B}(\is),\;t\geq0.
\end{equation*}
Of course ${U^*}^t\cdot a\otimes\bbbone_\mathcal{K}\cdot U^t=j_t[a]$ for every $a\in\mathcal{B}(\is)$
and $t\geq0$, as $\Theta$ commutes with $a\otimes\bbbone_\mathcal{K}$.

Given the state space $E$, we say that a QDS $T^t$ on some $\mathcal{B}(\is)$
extends a CMS $P^t$ on $\mathcal{L}^\infty(E)$ if 
\begin{itemize}
\item there is a $*$-isomorphism $f\mapsto m_f$ between 
$\mathcal{L}^\infty(E)$ and a (commutative) subalgebra of $\mathcal{B}(\is)$ such that
\begin{equation*}
    T[m_f] = m_{Pf}, \qquad \forall f \in \mathcal{L}^\infty(E).
\end{equation*}
\item for every $k\in E$ there exists a state $\rho_k\in\mathcal{S}(\is)$ such that 
\begin{equation}\label{eqexp}
    P^t\,f(k) = \tr\Big[T^t[m_f]\,\rho_k\Big], \qquad \forall f \in \mathcal{L}^\infty(E),
    \; t\geq0.
\end{equation}
\end{itemize}
Such extension always exists. For example, if we take $\is=\mathbb{C}^N=\mathbb{C}^{|E|}$ 
with its canonical basis $\{|i\rangle\}_{i\in E}$ and we
embed $\mathcal{L}^\infty(E)$ in $\mathcal{B}(\is)$ by
\begin{equation}\label{emb}
    f \mapsto m_f=\sum_{i\in E} f(i) \, |i\rangle\langle i|,
\end{equation}
$*$-isomorphism between $\mathcal{L}^\infty(E)$ and the subalgebra of the diagonal operators
$\mathcal{D}(\is)\subseteq\mathcal{B}(\is)$, then, using decomposition \eqref{dcc} and notations 
\eqref{detmatrix}, a stochastic matrix $P$ is extended by
\begin{equation}\label{Pext}
    T[a] = \sum_{\begin{subarray}{c} \scriptscriptstyle\ell\in L \\
    \scriptscriptstyle i\in E \end{subarray}}
    q_\ell\,|i\rangle\langle\beta_\ell(i)|\,a\,|\beta_\ell(i)\rangle\langle i|, \qquad
    \forall a\in\mathcal{B}(\is),
\end{equation}
and equation \eqref{eqexp} holds if every system state $k\in E$ is associated to the pure state
$\rho_k=|k\rangle\langle k|$.
Note that this extension maps $\mathcal{B}(\is)$ to $\mathcal{D}(\is)$ in only one step. 
The extension of a CMS is not unique at all: for example, with the same embedding \eqref{emb} and the
same choice of the states $\rho_k$, 
every permutation $P$ can be extended also by a $*$-automorphism
$T[a]=u^*\,a\,u$, provided that, for every $j$ in $E$, $u\,|j\rangle=P^*\,|j\rangle$ up to a phase 
factor (necessary and sufficient condition).

Given a CMS $P^t$ on $\mathcal{L}^\infty(E)$ and a quantum extension $T^t$ on some $\mathcal{B}(\is)$,
we say that a dilation $(\mathcal{K},U,\sigma)$ of $T^t$ extends a standard dilation 
$(G^\mathbb{Z},\mathcal{G},\alpha,Q_P)$ of $P^t$ if
\begin{itemize}
\item there exists a $*$-isomorphism $F\mapsto m_F$ between
$\LZ=\mathcal{L}^\infty(E)\otimes\Big(\bigotimes_{n\in\mathbb{Z}}\mathcal{L}^\infty(G_n)\Big)$ and a
(commutative) subalgebra of $\mathcal{B}(\is\otimes\mathcal{K})$ such that
\begin{equation*}
    {U^*}^t\,m_F\,U^t = m_{F\circ\alpha^t}, \qquad \forall F\in\LZ,\;t\in\mathbb{Z}
\end{equation*}
\item for every $k\in E$ there exists a state $\rho_k\in\mathcal{S}(\is)$ such that 
\begin{equation}\label{eqdistr}
    \mathbb{E}_{k,P}\big[\eta(F_{1},\ldots,F_{n})\big] =  
    \tr\big[\eta(m_{F_1},\ldots,m_{F_n})\,\rho_k\otimes\sigma\big],
\end{equation}
for every $F_1,\ldots,F_n \in\LZ$ and every 
$\eta:\mathbb{C}^n\to\mathbb{R}$ bounded and continuous.
\end{itemize} 
Note that the spatial tensor product of $C^*$-algebras $\LZ$ (\cite{KR} \textsection 11.4)
is a proper sub $C^*$-algebra of $\mathcal{L}^\infty(\mathcal{F})$ (it consists of all continuous
functions on the compact set $E\times G^\mathbb{Z}$) and that it is invariant for the evolution 
$F\mapsto F\circ\alpha^t$. Therefore also its image under $F\mapsto m_F$ has to be an abelian sub
$C^*$-algebra of $\mathcal{B}(\is\otimes\mathcal{K})$, invariant for the group of $*$-automorphisms 
$A\mapsto{U^*}^tAU^t$, which hence extends the classical evolution.
Moreover, the states $Q_P$ and $\sigma$ of the environments, together with the dynamics $\alpha$ and
$U$, give rise to the same joint distribution for the trajectories of the global random variables 
$F\in\LZ$ and for
the trajectories of the corresponding normal operators $m_F\in\mathcal{B}(\is\otimes\mathcal{K})$.
This happens for every possible starting state $k$ of the classical Markov system and for the
corresponding state $\rho_k$ of its quantum counterpart.

\begin{theorem}
Let $P$ be a stochastic matrix on a finite state space $E$ and let 
$(G^\mathbb{Z},\mathcal{G},\alpha,Q_P)$ be the dilation provided by Theorem \ref{cud} with
$Q_P=\bigotimes_\mathbb{Z}\widetilde{q}$. Then there exist a Hilbert space $\is$, a stochastic map
$T$ on $\mathcal{B}(\is)$ extending $P$, and a dilation $(\mathfrak{Z},V,\phi)$ of $T$ such that the
associated dilation $(\bigotimes_\mathbb{Z}\mathfrak{Z},U,\bigotimes_{\mathbb{Z}}\phi)$ of
$T^t$ extends $(G^\mathbb{Z},\mathcal{G},\alpha,\bigotimes_\mathbb{Z}\widetilde{q})$.
\end{theorem}

\noindent {\sl Proof.}
As in Theorem \ref{cud}, let $E=\{1,\ldots,N\}$, let 
$D_\ell=\sum_{i\in E}|\beta_\ell(i)\rangle\langle i|$ be the deterministic matrices corresponding to
the deterministic evolutions $\beta_\ell:E\to E$, as in Eq.~\eqref{detmatrix}, labelled by $\ell\in
L$. Then fix a probability $q_\ell$ on $L$ such that $P=\sum_{\ell\in L}q_\ell D_\ell$, as in
Eq.~\eqref{dcc}, and set $G=E\times L$.

Now, taken $\is=\mathbb{C}^{|E|}$ with its basis $\{|i\rangle\}_{i\in E}$, embedded
$\mathcal{L}^\infty(E)$ in $\mathcal{D}(\is)\subseteq\mathcal{B}(\is)$ by \eqref{emb},
let $P$ be extended by the stochastic map $T$ \eqref{Pext}.
Taken $\mathfrak{Z}=\mathbb{C}^{|G|}=\mathbb{C}^{|E\times L|}$ with its basis
$\{|j,\ell\rangle\}_{j\in E,\ell\in L}$, let $V$ be the unitary operator on
$\is\otimes\mathfrak{Z}$ which equals the adjoint of the permutation matrix associated
to the coupling $\varphi$ \eqref{coupling}, with respect to the basis 
$\{|i,j,\ell\rangle\}_{i,j\in E,\ell\in L}$:
\begin{equation*}
    V=\sum_{i,j\in E,\ell\in L} |\varphi(i,(j,\ell))\rangle\langle i,j,\ell|.
\end{equation*}
Taken the pure state $\phi=\sum_{\ell\in L}\sqrt{q_\ell}\,|1,\ell\rangle$ in $\mathfrak{Z}$, we are
going to show that
$(\mathfrak{Z},V,\phi)$ is a dilation of $T$ enjoying the required properties. 

For every $a\in\mathcal{B}(\is)$ and every
$\rho\in\mathcal{S}(\is)$, $\rho=\sum_{i,i'\in E}\rho_{ii'}\,|i\rangle\langle i'|$,
\begin{multline*}
    \tr\Big[V^*\cdot a\otimes\bbbone_\mathfrak{Z}\cdot
    V\cdot\rho\otimes|\phi\rangle\langle\phi|\Big] 
    = \sum_{\begin{subarray}{c} \scriptscriptstyle i,i'\in E \\
    \scriptscriptstyle \ell,\ell'\in L \end{subarray}}
    \sqrt{q_\ell}\,\sqrt{q_{\ell'}}\,\rho_{ii'}\,
    \langle i',1,\ell'| V^*\cdot a\otimes\bbbone_\mathfrak{Z}\cdot V |i,1,\ell\rangle \\
    = \sum_{\begin{subarray}{c} \scriptscriptstyle i,i'\in E \\
    \scriptscriptstyle \ell,\ell'\in L \end{subarray}}
    \sqrt{q_\ell}\,\sqrt{q_{\ell'}}\,\rho_{ii'}\,
    \langle \varphi(i',(1,\ell'))| a\otimes\bbbone_\mathfrak{Z}|\varphi(i,(1,\ell))\rangle \\
    = \sum_{\begin{subarray}{c} \scriptscriptstyle i,i'\in E \\
    \scriptscriptstyle \ell,\ell'\in L \end{subarray}}
    \sqrt{q_\ell}\,\sqrt{q_{\ell'}}\,\rho_{ii'}\,
    \langle \beta_{\ell'}(i'),1,\ell'| a\otimes\bbbone_\mathfrak{Z}|\beta_\ell(i),1,\ell\rangle \\
    = \sum_{i\in E, \ell\in L} q_\ell\,\rho_{ii}\,
    \langle \beta_{\ell}(i)|a|\beta_\ell(i)\rangle 
    =\tr\Big[T[a]\,\rho\Big],  
\end{multline*}
so that $(\mathfrak{Z},V,\phi)$ dilates $T$. 

Let now $(\mathcal{K},U,\Psi)$ be the associated dilation of $T^t$ given by
$\mathcal{K}=\bigotimes_{n\in\mathbb{Z}}\mathfrak{Z}_n$ with respect to the stabilizing sequence
$\psi_n\equiv\phi$, by the unitary operator $U$ defined by Eq.~\eqref{qshift}
and \eqref{glev} and by the pure state $\Psi=\bigotimes_{n\in\mathbb{Z}}\psi_n$. 
In order to show that
$(\mathcal{K},U,\Psi)$ extends $(G^\mathbb{Z},\mathcal{G},\alpha,\bigotimes_\mathbb{Z}\widetilde{q})$,
let $\mathcal{D}(\mathfrak{Z}_n)$ denotes the subalgebra of $\mathcal{B}(\mathfrak{Z}_n)$ 
consisting of
diagonal operators with respect to the basis 
$\{|j_n,\ell_n\rangle\}_{j_n\in E, \ell_n\in L}$. Then let $F \mapsto m_F$ be the $*$-isomorphism
between $\LZ$
and the spatial tensor product
$\DZ=\mathcal{D}(\is)\otimes\Big(\bigotimes_{n\in\mathbb{Z}}\mathcal{D}(\mathfrak{Z}_n)\Big)$
defined by mapping each
\begin{equation*}
    F\in \Lm=\mathcal{L}^\infty(E)\otimes\Big(\bigotimes_{n=-m}^m\mathcal{L}^\infty(G_n)\Big) = 
    \mathcal{L}^\infty\Big(E\times\big(\bigtimes_{n=-m}^m G_n\big)\Big), \qquad m\geq1,
\end{equation*} 
to the operator
\begin{equation*}
    m_F= \sum_{\begin{subarray}{c} \scriptscriptstyle i,j_{-m},\ldots,j_m\in E \\
    \scriptscriptstyle \ell_{-m},\ldots,\ell_m\in L\end{subarray}}
    F(i;j_{-m},\ell_{-m};\ldots;j_m,\ell_m)
    \, |i;j_{-m},\ell_{-m};\ldots;j_m,\ell_m\rangle\langle i;j_{-m},\ell_{-m};\ldots;j_m,\ell_m|
\end{equation*}
belonging to $\mathcal{D}(\is)\otimes\Big(\bigotimes_{n=-m}^m\mathcal{D}(\mathfrak{Z}_n)\Big)
\subseteq\mathcal{B}\Big(\is\otimes\big(\bigotimes_{n=-m}^m \mathfrak{Z}_n\big)\Big)$. 
For every $F$ belonging to an algebra $\Lm$,
\begin{equation*}
    \Theta^*\,m_F\,\Theta=m_{F\circ\vartheta}, \qquad
    V^*\,m_F\,V=m_{F\circ\varphi}, \qquad
    U^*\,m_F\,U=m_{F\circ\alpha},
\end{equation*}
so that, by continuity,
\begin{equation*}
    {U^*}^t\,m_F\,U^t=m_{F\circ\alpha^t}, \qquad\forall F\in\LZ,\; t\in\mathbb{Z}.
\end{equation*}

Finally let us consider the joint distribution of trajectories and let us show that property 
\eqref{eqdistr} holds if every system state $k\in E$ is associated to the pure state
$\rho_k=|k\rangle\langle k|$ on $\mathcal{B}(\is)$. Since $\eta(F_{1},\ldots,F_{n})$ belongs to $\LZ$
whenever $F_1,\ldots,F_n$ all belong to $\LZ$ and $\eta$ is bounded and continuous, it is to be checked
only that
\begin{equation*}
    \mathbb{E}_{k,P}[F] =  
    \tr\Big[m_{F}\cdot|k\rangle\langle k|\otimes|\Psi\rangle\langle\Psi|\Big], 
    \qquad \forall F\in\LZ.
\end{equation*}
If $F$ belongs to $\Lm$, then
\begin{multline*}
    \tr\Big[m_{F}\cdot|k\rangle\langle k|\otimes|\Psi\rangle\langle\Psi|\Big]
    =\sum_{\begin{subarray}{c} \scriptscriptstyle j_{-m},\ldots,j_m\in E \\
    \scriptscriptstyle \ell_{-m},\ldots,\ell_m\in L \end{subarray}}
    F(k;j_{-m},\ell_{-m};\ldots;j_m,\ell_m)\,\prod_{n=-m}^m
    \Big|\langle j_n,\ell_n|\phi\rangle\Big|^2 \\
    =\sum_{\begin{subarray}{c} \scriptscriptstyle j_{-m},\ldots,j_m\in E \\
    \scriptscriptstyle \ell_{-m},\ldots,\ell_m\in L \end{subarray}}
    F(k;j_{-m},\ell_{-m};\ldots;j_m,\ell_m)\,\prod_{n=-m}^m \widetilde{q}_{(j_n,\ell_n)}
    =\mathbb{E}_{k,P}[F],
\end{multline*}
and the general case follows by continuity.
\QED


\begin{thebibliography}{}
\bibitem{BMP}
Baldi, P.; Mazliak, L.; Priouret, P. \textsl{Martingales and Markov chains}. 
Chapman \& Hall/CRC, Boca Raton, FL, 2002.
\bibitem{C}
Chung, K.L.: \textsl{Markov chains with stationary transition probabilities}. 
Die Grundlehren der mathematischen Wissenschaften, Bd. \textbf{104} Springer-Verlag, 
Berlin-G\"ottingen-Heidelberg, 1960.
\bibitem{D61}
Davis, A.S. Markov chains as random input automata.  Amer. Math. Monthly  \textbf{68}  (1961) 
264--267.
\bibitem{G}
Ghom, R.: \textsl{Noncommutative Stationary Processes}. L.N.M. \textbf{1839} Berlin: Springer, 2004
\bibitem{KR}
Kadison, R.V.; Ringrose, J.R. \textsl{Fundamentals of the theory of
operator algebras. Vol. II. Advanced theory}. Pure and Applied Mathematics, 100. Academic Press, Inc.,
Orlando, FL, 1986.
\bibitem{K82}
K\"ummerer, B. A Dilation theory for completely positive operators on $W^*$-algebras. Thesis,
T\"ubingen, 1982.
\bibitem{K88}
K\"mmerer, B. Survey on a theory of noncommutative stationary Markov processes.
\textsl{Quantum probability and applications, III (Oberwolfach, 1987)}, 154--182,
Lecture Notes in Math., \textbf{1303}, Springer, Berlin, 1988. 
\bibitem{LM84}
Lewis, J.T.; Maassen, H. Hamiltonian models of classical and quantum stochastic processes.  
\textsl{Quantum probability and applications to the quantum theory of irreversible processes 
(Villa Mondragone, 1982)},  245--276, Lecture Notes in Math., \textbf{1055}, Springer, Berlin, 1984.
\bibitem{M}
Meyer, P.A.: \textsl{Quantum probability for probabilists}.
Lecture Notes in Mathematics, \textbf{1538}.
Springer-Verlag, Berlin, 1993.
\bibitem{P}
Parthasarathy, K.R.: \textsl{An Introduction to Quantum
Stochastic Calculus}. Basel-Boston-Berlin: Birkh\"{a}user, 1992.
\bibitem{P89} 
Parthasarathy, K.R. Discrete time quantum stochastic flows, Markov chains and chaos expansions. 
Probability theory and mathematical statistics, Vol. II (Vilnius, 1989),  288--297, "Mokslas",
Vilnius, 1990.
\end{thebibliography}
\end{document}